\documentclass[draft]{amsart}
\usepackage{amssymb,graphicx,enumerate,amsthm}
\usepackage{multirow}
\usepackage{xcolor}
\usepackage{cite}

\newcommand{\Z}{\mathbb{Z}}

\numberwithin{equation}{section}
\newtheorem{theorem}{Theorem}[section]

\newtheorem{conjecture}[theorem]{Conjecture}

\begin{document}

\makeatletter
\def\imod#1{\allowbreak\mkern10mu({\operator@font mod}\,\,#1)}
\makeatother

\author{Alexander Berkovich}
   	\address{Department of Mathematics, University of Florida, 358 Little Hall, Gainesville FL 32611, USA}
   	\email{alexb@ufl.edu}

\thanks{Research is partly supported by the Simons Foundation, Award ID: 308929.} 
%Research of the second author is supported by the Austrian Science Fund FWF, SFB50-07 and SFB50-09 Projects.}

\title[\scalebox{.9}{Some New Positive Observations}]{Some New Positive Observations}
%New positivity preserving transformations for \lowercase{$q$}-binomial coefficients
     
\begin{abstract}
We revisit Bressoud's generalized Borwein conjecture. Making use of the new positivity-preserving transformations for 
$q$-binomial coefficients we establish the truth of infinitely many cases of the Bressoud conjecture.
In addition, we prove new bounded versions of Lebesgue's identity and of Euler's Pentagonal Number Theorem.
Finally, we discuss new companions to Andrews-Gordon $mod\text{ }21$ and Bressoud $mod\text{ }20$ identities.
\end{abstract}

  

%11C08		Polynomials
%11P81  	Elementary theory of partitions [See also 05A17]
%11P82  	Analytic theory of partitions
%11P83  	Partitions; congruences and congruential restrictions
%11P84  	Partition identities; identities of Rogers-Ramanujan type
%%11B65  	Binomial coefficients; factorials; $q$-identities
%05A10  	Factorials, binomial coefficients, combinatorial functions [See also 11B65, 33Cxx]
%05A15  	Exact enumeration problems, generating functions [See also 33Cxx, 33Dxx]
%05A17  	Partitions of integers [See also 11P81, 11P82, 11P83]

\date{\today}
   
\dedicatory{Dedicated to the memory of Dora Bitman}
   
\maketitle

\section{Introduction}

Bressoud \cite{BR3} considered the following polynomials

\begin{equation}
G(N,M,\alpha,\beta,K,q)=\sum_{j=-\infty}^\infty(-1)^j q^{Kj\frac{(\alpha+\beta)j+(\alpha-\beta)}{2}}{N+M \brack N-Kj}_q,
\label{1.1}
\end{equation}
where
\begin{equation}
{m+n \brack m}_q:= \left\lbrace \begin{array}{ll}\frac{(q)_{m+n}}{(q)_m(q)_{n}},&\text{for }m,n\in\mathbf N,\\
   0,&\text{otherwise,}\end{array}\right.
\label{1.2-binom_def}
\end{equation}
and
\begin{equation*}
(q)_m=\prod_{j=1}^m(1-q^j),\text{for }m\in\mathbf N,
\end{equation*}
where $\mathbf N$ denotes the set of nonnegative integers.
More generally, for $m\in\mathbf N$ we define
\begin{align}
(a)_m=(a;q)_m &=\prod_{j=0}^{m-1}(1-aq^j),\nonumber\\
(a_1,a_2,\ldots,a_k;q)_m &=(a_1;q)_m(a_2;q)_m\ldots(a_k;q)_m.
\label{1.3}
\end{align}
Here and throughout we assume that $|q|<1$. We note that $(a)_0=1$.\\

In 1996, Bressoud \cite{BR3} conjectured that
\begin{conjecture}\label{Bressoud conjecture}
Let $K\in\Z_{>1}$, $N,M,\alpha K,\beta K\in\Z_{\geq0}$ such that 
\begin{align}
1 & \leq\alpha + \beta\leq 2K-1,\nonumber \\
\beta-K & \leq N-M\leq K-\alpha,
\label{1.4}
\end{align}
(strict inequality when $K=2$). Then $G(N,M,\alpha,\beta,K,q)\geq0$.
\end{conjecture}
\noindent
Here, and everywhere, $P(q)\geq0$ means that $P(q)$ is a polynomial in $q$ with nonnegative coefficients.\\
We remark that
\begin{equation*}
{m+n \brack m}_q\geq 0.
\end{equation*}
The famous conjecture of Peter Borwein (Theorem since 2019 \cite{CW}) can be stated as
\begin{align}
A_n(q) & = G(n,n,\frac{5}{3},\frac{4}{3},3,q)\geq 0,\nonumber\\
B_n(q) & = G(n-1,n+1,\frac{7}{3},\frac{2}{3},3,q)\geq 0,\nonumber \\
C_n(q) & = G(n-1,n+1,\frac{8}{3},\frac{1}{3},3,q)\geq 0,
\label{1.5}
\end{align}
and
\begin{equation}
\prod_{k=1}^n(1-q^{3k-1})(1-q^{3k-2})=A_n(q^3)-qB_n(q^3)-q^2C_n(q^3).
\label{1.6}
\end{equation}
When $\alpha,\beta\in\Z$, $G(N,M,\alpha,\beta,K,q)$ is a generating function for the so-called partitions with prescribed hook differences \cite{ABBBFV}.
Bressoud's conjecture is nontrivial when $\alpha,\beta$ assume fractional values.\\
Many cases of Bressoud's conjecture were settled in the literature \cite{BR2},\cite{IKS},\cite{W1},\cite{W2},\cite{BW},\cite{CW}.\\
In the next section, we will show how to settle new infinite family of cases. 
\begin{theorem}\label{theorem1}
For $L\in \mathbf N$, $\nu\in\Z_{>0}$, $s=0,1,2,\ldots,\nu-1$
\begin{equation}
G(L,L+1+2s,(\nu+1)(1+\frac{1+2s}{2\nu+1}),(\nu+1)(1-\frac{1+2s}{2\nu+1}),2\nu+1,q)\geq0.
\label{1.7}
\end{equation}
\end{theorem}
Also, in Section 2, we discuss new bounded versions of Lebesgue's identity and of Euler's Pentagonal Number Theorem.
In Section 3, we establish and prove some additional isolated positivity results and introduce new companions to Andrews-Gordon $mod\text{ }21$ 
and Bressoud $mod\text{ }20$ identities.\\
We conclude this section with a list of seven useful formulas, which can be found in \cite{A1}:
\begin{align}
\label{1.8-Binom_limit}
\lim_{L\rightarrow\infty}{L\brack m}_q &= \frac{1}{(q)_m},\\
\label{1.9-Binom_limit2} \lim_{L,M\rightarrow\infty}{L+M \brack L}_q &= \frac{1}{(q)_\infty},
\end{align} 
\begin{equation}
{n+m \brack n}_{q^{-1}}=q^{-nm}{n+m \brack n}_q,
\label{1.10} 
\end{equation}
\begin{equation}
{n \brack m}_q={n-1 \brack m-1}_q+q^m{n-1 \brack m}_q = {n-1 \brack m}_q+q^{n-m}{n-1 \brack m-1}_q,
\label{1.11} 
\end{equation}
\begin{equation}
\sum_{n\geq 0} q^{n \choose 2}z^n{L \brack n}_q = (-z;q)_L,
\label{1.12} 
\end{equation}
\begin{equation}
\sum_{j=-\infty}^\infty (-1)^j z^j q^{j^2} = \left(q^2,\frac{q}{z},zq;q^2\right)_\infty,
\label{JTP}
\end{equation} 
\begin{equation}
\sum_{j=-\infty}^\infty (-1)^j q^{j^2} z^j {L+M \brack L-j}_{q^2}  = \left(\frac{q}{z};q^2\right)_M \left(zq;q^2\right)_L,
\label{1.14}
\end{equation} 
with $L,M,m,n\in\mathbf N$.

\section{Positivity-preserving Transformations}

We start with the following summation formula
\begin{theorem}
For $L\in\mathbf N$, $a\in\mathbf Z$
\begin{equation}
\sum_{k\geq0} C_{L,k}(q){k \brack \lfloor\frac{k-a}{2}\rfloor}_q=q^{T(a)}{2L+1 \brack L-a}_q,
\label{2.1}
\end{equation}
\end{theorem}
\noindent
where 
\begin{equation*}
T(j):={j+1\choose 2}
\end{equation*}
and
\begin{equation}
C_{L,k}(q)=\sum_{m=0}^L q^{T(m)+T(m+k)}{L \brack m,k}_q,      
\label{2.2}
\end{equation}
with 
\begin{equation}
{L \brack m,k}_q={L \brack m}_q {L-m \brack k}_q={L \brack k}_q {L-k \brack m}_q\geq 0. 
\label{2.3}
\end{equation}
Observe that $C_{L,k}(q)\geq0$.
Using transformation \eqref{2.1} it is easy to check that identity of the form
\begin{equation}
F(L,q)=\sum_{j=-\infty}^\infty\alpha(j,q){L \brack \lfloor\frac{L-j}{2}\rfloor}_q,
\label{2.4}
\end{equation}
implies that the following identity holds
\begin{equation}
\sum_{k\geq0}C_{L,k}(q)F(k,q)=  
\sum_{j=-\infty}^\infty\alpha(j,q)\sum_{k\geq0}C_{L,k}(q){k \brack \lfloor\frac{k-j}{2}\rfloor}_q=
\sum_{j=-\infty}^\infty\alpha(j,q)q^{T(j)}{2L+1 \brack L-j}_q.
\label{2.5}
\end{equation}
Hence, if $F(L,q)\geq0$ then
\begin{equation}
\sum_{j=-\infty}^\infty\alpha(j,q)q^{T(j)}{2L+1 \brack L-j}_q\geq0.
\label{2.6}
\end{equation}
For that reason, we say that \eqref{2.1} is positivity-preserving. \\
Transformation \eqref{2.1} is an easy corollary of the theorem proven in \cite{BU}. 
\begin{theorem}[Berkovich--Uncu]
\begin{equation}
\sum_{k\geq 0} q^{T(k)} {L\brack k}_q \left\{ T_{-1}\left(\begin{array}{c}k\\ a \end{array};q \right) + T_{-1}\left(\begin{array}{c}k\\ a+1 \end{array};q \right) \right\} = q^{T(a)} {2L+1\brack L-a}_q.
\label{2.6a}
\end{equation}
\end{theorem}
\noindent
The Andrews-Baxter $q$-trinomial coefficients \cite{AB} can be defined as 
\begin{equation}
T_{-1} \left(\begin{array}{c}k\\ a \end{array};q \right) = 
\sum_{\substack{m\geq0,\\m\equiv k+a(mod2)}} q^{T(m)} {k\brack m}_q {k-m\brack \frac{k-m-a}{2}}_q.
\label{2.6b}
\end{equation}
It is easy to check that
\begin{equation}
 T_{-1}\left(\begin{array}{c}k\\ a \end{array};q \right) + T_{-1}\left(\begin{array}{c}k\\ a+1 \end{array};q \right) = 
\sum_{m\geq0} q^{T(m)} {k\brack m}_q {k-m\brack \lfloor\frac{ k-m-a}{2}\rfloor}_q.
\label{2.6c}
\end{equation}
Substituting \eqref{2.6c} into left hand side of \eqref{2.6a} and changing $k\rightarrow k+m$ we complete the proof of \eqref{2.1}.\\

It is instructive to compare \eqref{2.1} with the Corollary (2.6) in \cite{W2}.
\begin{theorem}[Warnaar]\label{Warnaar}
For $L\in\mathbf N$, $a\in\mathbf Z$
\begin{equation}
\sum_{k\geq0} W_{L,k}(q){2k \brack k-a}_q=q^{2a^2}{2L \brack L-2a}_q,
\label{2.6d}
\end{equation}
where
\begin{equation}
W_{L,k}(q)=\sum_{m=0}^L q^{(m+k)^2+k^2}{L \brack m,2k}_q\geq0.      
\label{2.6e}
\end{equation}
\end{theorem}
\noindent
Observe that unlike \eqref{2.6d}, transformation \eqref{2.1} can not be iterated. Interestingly enough, there exists an \textit{odd} 
companion to Theorem~\ref{Warnaar}.
\begin{theorem}\label{Odd_companion}
For $L\in\mathbf N$, $a\in\mathbf Z$
\begin{equation}
\sum_{k\geq0} O_{L,k}(q){2k+1 \brack k-a}_q=q^{4T(a)}{2L \brack L-2a-1}_q,
\label{2.6f}
\end{equation}
where
\begin{equation}
O_{L,k}(q)=\sum_{m=0}^L q^{2T(m+k)+2T(k)}{L \brack m,2k+1}_q\geq0.      
\label{2.6g}
\end{equation}
\end{theorem}
\noindent
We remark that while Theorem~\ref{Odd_companion} is not explicitly stated in \cite{W2}, it is a special case of an identity on page 222 there.\\

Schur's bounded version of Euler's Pentagonal Number Theory states
\begin{equation}
1=\sum_{j=-\infty}^\infty (-1)^j q^{\frac{3j+1}{2}j} {L \brack \lfloor\frac{L-3j}{2}\rfloor}_q.
\label{2.13}
\end{equation}
With the aid of \eqref{2.5} we can convert \eqref{2.13} into
\begin{equation}
0\leq\sum_{k=0}^L C_{L,k}(q)=\sum_{j=-\infty}^\infty (-1)^j q^{2j(3j+1)} {2L+1 \brack L-3j}_q.
\label{2.14}
\end{equation}
Hence, 
\begin{equation*}
G(L,L+1,\frac{8}{3},\frac{4}{3},3,q)\geq0.
\end{equation*}
Making use of \eqref{1.12}, it is easy to check that
\begin{equation}
\sum_{k=0}^L C_{L,k}(q)=
\sum_{k=0}^L q^{T(k)}{L \brack k}_q(-q)_k.
\label{2.15}
\end{equation}
And so identity \eqref{2.14} can be rewritten as
\begin{equation}
\sum_{k=0}^L q^{T(k)}{L \brack k}_q(-q)_k=\sum_{j=-\infty}^\infty (-1)^j q^{2j(3j+1)} {2L+1 \brack L-3j}_q.
\label{2.16}
\end{equation}
Letting $L\rightarrow\infty$ and using the Jacobi triple product identity \eqref{JTP} yields a special case of the Lebesgue identity \cite{L}
\begin{equation}
\sum_{m\geq0}\frac{q^{T(m)}}{(q)_m}(-q)_m=\frac{(q^4;q^4)_\infty}{(q)_\infty},
\label{2.17}
\end{equation}
and so, \eqref{2.16} is a new bounded version of the Lebesgue identity.
Perform $q\rightarrow\frac{1}{q}$ in \eqref{2.16} and use \eqref{1.10} together with 
\begin{equation}
(-q^{-1};q^{-1})_n=(-q)_n q^{-T(n)},\text{  }n\in\mathbf N
\label{2.18}
\end{equation}
to obtain after simplification a new polynomial version of Euler's Pentagonal Number Theorem
\begin{equation}
\sum_{k=0}^L (-q)_{L-k}q^{(L+1)k}{L \brack k}_q=\sum_{j=-\infty}^\infty (-1)^j q^{3j^2+j} {2L+1 \brack L-3j}_q.
\label{2.19}
\end{equation}
It proves that
\begin{equation}
G(L,L+1,\frac{4}{3},\frac{2}{3},3,q)\geq0.
\label{2.20}
\end{equation}

We now move on to prove Theorem~\ref{theorem1}. We start with the finite analogue of the Andrews-Gordon identity due to Foda-Quano \cite{FQ}. \\

For $L\in \mathbf N$, $\nu\in\Z_{>0}$, $s=0,1,\ldots,\nu-1$ 
\begin{align}
& \sum_{n_2,\ldots,n_\nu\geq0}q^{N_2^2+\ldots+N_{\nu}^2+N_{\nu+1-s}+\ldots+N_\nu}
\prod_{i=2}^\nu{n_i+L-2\sum_{j=2}^i N_j-E_{i,s}^\nu\brack n_i}_q= \nonumber\\
& \sum_{j=-\infty}^\infty (-1)^j q^{\frac{(2\nu+1)j^2+j(1+2s)}{2}} {L \brack \lfloor\frac{L-(2\nu+1)j-s}{2}\rfloor}_q.
\label{2.21}
\end{align}
Here, $N_j=n_j+n_{j+1}+\ldots+n_\nu$, $j=2,\ldots,\nu$ and $E_{i,s}^\nu=max(i+s-\nu,0)$. 
Observe that \eqref{2.13} is the case $\nu=1$ of \eqref{2.21}. \\

With the aid of \eqref{2.5} we obtain
\begin{align}
0\leq & \sum_{k,n_2,\ldots,n_\nu\geq0} C_{L,k}(q)q^{N_2^2+\ldots+N_{\nu}^2+N_{\nu+1-s}+\ldots+N_\nu} 
\prod_{i=2}^\nu{n_i+k-2\sum_{j=2}^i N_j-E_{i,s}^\nu\brack n_i}_q=\nonumber \\
& q^{T(s)}\sum_{j=-\infty}^\infty (-1)^j q^{(\nu+1)(2\nu+1)j^2+(\nu+1)(2s+1)j} {2L+1 \brack L-s-(2\nu+1)j}_q.
\label{2.22}
\end{align}
Hence, 
\begin{equation*}
G(L,L+1+2s,(\nu+1)(1+\frac{1+2s}{2\nu+1}),(\nu+1)(1-\frac{1+2s}{2\nu+1}),2\nu+1,q)\geq0,
\end{equation*}
for all $L\in\mathbf N$, $\nu\in\Z_{>0}$, $s=0,1,2,\ldots,\nu-1$. This completes the proof of Theorem~\ref{theorem1}.

\section{Further Observations}

We replace $q^2$ by $q$ in \eqref{1.14} and then
set $M=L,L+1$, $z=q^{\frac{1}{2}}$ to find that for $L\in\mathbf N$
\begin{equation}
\sum_{j=-\infty}^\infty (-1)^j q^{T(j)}{2L \brack L-j}_q=\delta_{L,0}
\label{3.1}
\end{equation}
and
\begin{equation}
\sum_{j=-\infty}^\infty (-1)^j q^{T(j)}{2L+1 \brack L-j}_q=0,
\label{3.2}
\end{equation}
where
$\delta_{L,0}=1$ if $L=0$ and $\delta_{L,0}=0$ if $L>0$.
The formulas \eqref{3.1} and \eqref{3.2} can be combined into 
\begin{equation}
\sum_{j=-\infty}^\infty (-1)^j q^{T(j)}{L \brack \lfloor\frac{L-2j}{2}\rfloor}_q=\delta_{L,0}.
\label{3.3}
\end{equation}
Applying Theorem~\ref{Warnaar} to \eqref{3.1} yields
\begin{equation}
W_{L,0}(q)=\sum_{n\geq0} q^{n^2}{L \brack n}_q=\sum_{j=-\infty}^\infty (-1)^j q^{\frac{5j+1}{2}j}{2L \brack L-2j}_q,
\label{3.4}
\end{equation}
which is Bressoud's bounded version of the first Rogers-Ramanujan identity \cite{BR2}.
Analogously, applying \eqref{2.5} to \eqref{3.3} yields
\begin{equation}
C_{L,0}(q)=\sum_{n\geq0} q^{n^2+n}{L \brack n}_q=\sum_{j=-\infty}^\infty (-1)^j q^{\frac{5j^2+3j}{2}}{2L+1 \brack L-2j}_q,
\label{3.5}
\end{equation}
which can be recognized as Warnaar's bounded version of the second Rogers-Ramanujan identity \cite{W1}.
Next, we perform the change of summation variables below 
\begin{align}
&\sum_{j=-\infty}^\infty (-1)^j q^{5T(j)}{2L \brack L-2j-1}_q=\nonumber\\
&-\sum_{j=-\infty}^\infty (-1)^j q^{5T(-1-j)}{2L \brack L+2j+1}_q=\\
&-\sum_{j=-\infty}^\infty (-1)^j q^{5T(j)}{2L \brack L-2j-1}_q\nonumber
\label{3.6}
\end{align}
to conclude that 
\begin{equation}
q^{L+1}\sum_{j=-\infty}^\infty (-1)^j q^{5T(j)}{2L \brack L-2j-1}_q=0.
\label{3.7}
\end{equation}
Adding \eqref{3.4} and \eqref{3.7} and employing recursion relation \eqref{1.11} we obtain
\begin{equation}
\sum_{n\geq0} q^{n^2}{L \brack n}_q=\sum_{j=-\infty}^\infty (-1)^j q^{\frac{5j^2+j}{2}}{2L+1 \brack L-2j}_q.
\label{3.8}
\end{equation}
Observe that \eqref{3.4} and \eqref{3.8} imply that for $k\in \mathbf N$
\begin{equation}
\sum_{n\geq0} q^{n^2}{\lfloor\frac{k}{2}\rfloor \brack n}_q=
\sum_{j=-\infty}^\infty (-1)^j q^{\frac{5j^2+j}{2}}{k \brack \lfloor\frac{k-4j}{2}\rfloor}_q.
\label{3.9}
\end{equation}
Apply Theorem~\ref{2.1} to \eqref{3.9} to obtain
\begin{equation}
\sum_{k,n\geq0}C_{L,k} q^{n^2}{\lfloor\frac{k}{2}\rfloor \brack n}_q=
\sum_{j=-\infty}^\infty (-1)^j q^{\frac{21j^2+5j}{2}}{2L+1 \brack L-4j}_q,
\label{3.10}
\end{equation}
which proves that 
\begin{equation*}
G(L,L+1,\frac{13}{4},2,4,q)\geq0.
\end{equation*}
In the limit as $L\rightarrow\infty$ \eqref{3.10} becomes
\begin{equation}
\sum_{m,k,n\geq0} \frac{q^{T(m)+T(m+k)+n^2} {\lfloor\frac{k}{2}\rfloor \brack n}_q}{(q)_m(q)_k}=
\frac{(q^{21},q^8,q^{13};q^{21})_n}{(q)_\infty}.
\label{3.11}
\end{equation}
This is to be contrasted with Andrews-Gordon identity mod $21$ \cite{A}
\begin{equation}
\sum_{n_1,n_2,\ldots,n_9\geq0}\frac{q^{N_1^2+N_2^2+\ldots+N_9^2+N_8+N_9}}{(q)_{n_1}(q)_{n_2}\ldots(q)_{n_9}}=
\frac{(q^{21},q^8,q^{13};q^{21})_\infty}{(q)_\infty},
\label{3.12}
\end{equation}
with $N_i=n_i+\ldots+n_9$, $i=1,\ldots,9$.
On the left of \eqref{3.11} one has $3$-fold sum, while on the left side of \eqref{3.12} one has $9$-fold sum.
Analogously, applying Theorem~\ref{Warnaar} to \eqref{3.4} and Theorem~\ref{Odd_companion} to \eqref{3.5} and \eqref{3.8},
we prove that
\begin{align*}
G(L,L,\frac{11}{4},\frac{5}{2},4,q)\geq 0,\\
G(L-1,L+1,4,\frac{5}{4},4,q)\geq 0, \\
G(L-1,L+1,\frac{15}{4},\frac{3}{2},4,q)\geq 0,
\end{align*}
and obtain, as $L\rightarrow\infty$
\begin{equation}
\sum_{m,k,n\geq0} \frac{q^{k^2+(m+k)^2+n^2}}{(q)_m(q)_{2k}} {k \brack n}_q=
\frac{(q^{21},q^{10},q^{11};q^{21})_\infty}{(q)_\infty},
\label{3.13}
\end{equation}
\begin{equation}
\sum_{m,k,n\geq0} \frac{q^{2T(k)+2T(m+n)+2T(n)}}{(q)_m(q)_{2k+1}} {k \brack n}_q=
\frac{(q^{21},q^5,q^{16};q^{21})_\infty}{(q)_\infty},
\label{3.14}
\end{equation}
and
\begin{equation}
\sum_{m,k,n\geq0}\frac{q^{2T(k)+2T(m+n)+n^2}}{(q)_m(q)_{2k+1}} {k \brack n}_q=
\frac{(q^{21},q^6,q^{15};q^{21})_\infty}{(q)_\infty},
\label{3.15}
\end{equation}
respectively.\\

In \cite[p. 2332]{BW} the following identity was derived
\begin{equation}
\sum_{j=-\infty}^\infty (-1)^j q^{2j^2}{2L \brack L-2j}_q=(-q;q^2)_L.
\label{3.16x}
\end{equation}
We now follow a well-trodden path and check that
\begin{equation}
q^{L+1}\sum (-1)^j q^{2j^2+2j}{2L \brack L-2j-1}_q=0.
\label{3.17y}
\end{equation}
Adding \eqref{3.16x} and \eqref{3.17y} we derive, with the aid of \eqref{1.11}, that
\begin{equation}
\sum_{j=-\infty}^\infty (-1)^j q^{2j^2}{2L+1 \brack L-2j}_q=(-q;q^2)_L.
\label{3.18z}
\end{equation}
Equations \eqref{3.16x} and \eqref{3.18z} imply that for $k\in\mathbf N$
\begin{equation}
\sum_{j=-\infty}^\infty (-1)^j q^{2j^2}{k \brack \lfloor\frac{k-4j}{2}\rfloor}_q=(-q;q^2)_{\lfloor\frac{k}{2}\rfloor}.
\label{3.19w}
\end{equation}
Applying Theorem~\ref{2.1} to \eqref{3.19w} and letting $L\rightarrow\infty$ we obtain
\begin{equation}
\sum_{m,k\geq0}\frac{q^{T(m)+T(m+k)}(-q;q^2)_{\lfloor\frac{k}{2}\rfloor}}{(q)_m(q)_k}=
\frac{(q^{20},q^8,q^{12};q^{20})_\infty}{(q)_\infty}.
\label{3.20}
\end{equation}
Compare it with the Bressoud formula in \cite{BR1}
\begin{equation}
\sum_{n_1,\ldots,n_9\geq0}\frac{q^{N_1^2+\ldots+N_9^2+N_8+N_9}}{(q)_{n_1}\ldots(q)_{n_8}(q^2;q^2)_{n_9}}=
\frac{(q^{20},q^8,q^{12};q^{20})_\infty}{(q)_\infty},
\label{3.21}
\end{equation}
where $N_i=n_i+\ldots+n_9$, $i=1,\ldots,9$.\\

Analogously, applying Theorem~\ref{Warnaar} to \eqref{3.16x} and Theorem~\ref{Odd_companion} to \eqref{3.18z} we get as $L\rightarrow\infty$
\begin{equation}
\sum_{m,k\geq0}\frac{q^{k^2+(m+k)^2}}{(q)_m(q)_{2k}}(-q;q^2)_k=\frac{(q^{20},q^{10},q^{10};q^{20})_\infty}{(q)_\infty}
\label{3.22}
\end{equation}
and
\begin{equation}
\sum_{m,k\geq0}\frac{q^{2T(k)+2T(m+k)}}{(q)_m(q)_{2k+1}}(-q;q^2)_k=\frac{(q^{20},q^6,q^{14};q^{20})_\infty}{(q)_\infty},
\label{3.23}
\end{equation}
respectively.\\
For our final example, we employ Dyson's identity \cite{BA}, \cite[p. 2330]{BW}
\begin{equation}
\sum_{j=-\infty}^\infty (-1)^j q^{T(3j)} {2L+1 \brack L-3j}_q = \frac{(q^3;q^3)_L}{(q)_L}.
\label{3.24}
\end{equation}
Applying Theorem~\ref{Odd_companion} to \eqref{3.24} yields
\begin{equation}
\sum_{j=-\infty}^\infty (-1)^j q^{5T(3j)} {2L \brack L-1-6j}_q = \sum_{k\geq0} O_{L,k} \frac{(q^3;q^3)_k}{(q)_k}\geq0.
\label{3.25}
\end{equation}
This proves that 
\begin{equation}
G(L-1,L+1,5,\frac{5}{2},6,q)\geq0.
\label{3.26}
\end{equation}
Letting $L\rightarrow\infty$ in \eqref{3.25} and using \eqref{JTP} we arrive at a new elegant result
\begin{equation}
\sum_{m,k\geq0}\frac{q^{2T(k)+2T(m+k)}}{(q)_m(q)_{2k+1}}\frac{(q^3;q^3)_k}{(q)_k}=\frac{(q^{15};q^{15})_\infty}{(q)_\infty}.
\label{3.27}
\end{equation}

\section{Acknowledgement}

I would like to thank George Andrews, James Mc Laughlin and Ali Uncu for their kind interest.


\begin{thebibliography}{99}

\bibitem{A}G. E. Andrews, \textit{An analytic generalization of the Rogers-Ramanujan identities for odd moduli}, 
           Proc. Nat. Acad. Sci. USA \textbf{71} (1974), 4082--4085.

\bibitem{A1}G. E. Andrews, \textit{The theory of partitions}, Cambridge Mathematical Library, Cambridge University Press, Cambridge, 1998            Reprint of the 1976 original. MR1634067 (99c:11126)

\bibitem{AB}G. E. Andrews, and R. J. Baxter, \textit{Lattice gas generalization of the hard hexagon model. III. q-Trinomial coefficients}, 
            J. Statist. Phys. \textbf{47} (1987), no: 3-4, 297--330.

\bibitem{ABBBFV}G. E. Andrews, R. J. Baxter, D. M. Bressoud, W. H. Burge, P. J. Forrester, G. Viennot, 
                \textit{Partitions with prescribed hook differences}, Europe J. Comb. \textbf{8} (1987) no. 4, 341--350. 

\bibitem{BA}W. N. Bailey, \textit{Some identities in combinatorial analysis}, Proc. London Math. Soc. \textbf{2} (1947), 421--435. 

\bibitem{BU}A. Berkovich, A. K. Uncu, \textit{Elementary polynomial identities involving $q$-trinomial coefficients}, 
            Ann. Comb. \textbf{23} (2019), 549--560. 

\bibitem{BW}A. Berkovich, S. O. Warnaar, \textit{Positivity preserving transformations for $q$-binomial coefficients}, 
            Trans. Amer. Math. Soc. \textbf{357} (2005), no. 6, 2291--2351.

\bibitem{BR1}D. M. Bressoud, \textit{An analytic generalization of Rogers-Ramanujan identities with interpretation},  
             Quart. J. Math. Oxford (2) \textbf{31} (1980), 389--399. 

\bibitem{BR2}D. M. Bressoud, \textit{Some identities for terminating $q$-series}, 
            Math. Proc. Camb. Phil. Soc. \textbf{89} (1981), 211--223.              

\bibitem{BR3}D. M. Bressoud, \textit{The Borwein conjecture and partitions with prescribed hook differences},  
            Electron. J. Combin. \textbf{3} (1996) no. 2, Research Paper 4. 

\bibitem{FQ}O. Foda, Y. H. Quano, \textit{Polynomial identities of the Rogers-Ramanujan type}, 
            Int. J. Mod. Phys. A \textbf{10} (1995), 2291--2315. 

\bibitem{IKS}M. E. H. Ismail, D. Kim, D. Stanton, \textit{Lattice paths and positive trigonometric sums}, 
             Const. Approx. \textbf{15} (1999), 69--81.							
	
\bibitem{L}V. A. Lebesgue, \textit{Sommation de quelques series}, 
           J. Math. Pure. Appl. \textbf{5} (1840), 42--71. 							
		
\bibitem{CW}Chen Wang, \textit{Analytic proof of the Borwein conjecture}, arXiv:1901.10886v1. [math.CO].							
								
\bibitem{W1}S. O. Warnaar, \textit{The generalized Borwein conjecture. I.The Burge transform, in: B.C. Berndt, K. Ono (Eds.), 
            $q$-Series with Applications to Combinatorics, Number Theory and Physics}, Contemp. Math. \textbf{291}, 
						AMS, Providence, RI (2001), 243--267. 
	
\bibitem{W2}S. O. Warnaar, \textit{The generalized Borwein conjecture. II. Refined $q$-trinomial coefficients}, 
            Discrete Math. \textbf{272} (2003), no. 2-3, 215--258.

\end{thebibliography}
\end{document}